\documentclass[12pt]{article}
\usepackage{amsfonts}
\usepackage{amsmath}
\usepackage{listings}
\usepackage{epsfig}
\usepackage{xcolor}
\usepackage{geometry}
\title{Vertex-Minimal Paper Tori}
\author{Richard Evan Schwartz \thanks{\hskip 5 pt Supported by 
N.S.F. Research Grant DMS-2505281}}

\newtheorem{theorem}{Theorem}[section]

\newtheorem{lemma}[theorem]{Lemma}

\newtheorem{corollary}[theorem]{Corollary}

\def\startproof{{\bf {\medskip}{\noindent}Proof: }}

\def\endproof{$\spadesuit$  \newline}

\def\H{\mbox{\boldmath{$H$}}}%
\def\R{\mbox{\boldmath{$R$}}}%

\begin{document}

\maketitle

\begin{abstract}
A {\it paper torus\/} is an embedded polyhedral
torus that is isometric to a flat torus in the intrinsic sense.
  We prove that there
  does not exist a paper torus with $7$ vertices,
  and that there does exist a paper torus with $8$ vertices.
\end{abstract}

\section{Introduction}

\subsection{Context and Main Results}

A {\it flat torus\/} is a quotient of the
form $\R^2/\Lambda$, where $\Lambda$ is a
lattice of translations of $\R^2$.  More concretely,
a flat torus is what you get when you identify the
opposite sides of a parallelogram by translations.

A {\it polyhedral torus\/} is a continuous, piecewise
affine embedding
$\phi: T \to \Omega \subset \R^3.$
Here $T$ is a flat torus that has been triangulated,
and $\phi$ is a continuous embedding that
is affine on each triangle of the triangulation.
If, additionally, $\phi$ is an affine isometry on each
triangle, we say that $\phi$ is {\it flat\/}.
When $\phi$ is a flat embedding, $\phi$ tells how to build
a torus in $\R^3$ out of finitely many triangles so that
the cone angle around each vertex
is $2 \pi$.  We call a flat embedded polyhedral torus
a {\it paper torus\/}.   Some authors would call this
an origami torus.

Surprisingly, paper tori exist.
The 1960 paper
of Y. Burago and V. Zalgaller 
[{\bf BZ1\/}] gives the first constructions.
Some of the simplest paper tori are
called {\it diplotori\/}.   These are
described by U. Brehm [{\bf Br\/}] in 1978,
and also described
in H. Segerman's book [{\bf Se\/}, \S 6].

The 1995 paper [{\bf BZ2\/}] proves that
one can realize every isometry class of flat torus as a paper torus.
The works of T. Tsuboi [{\bf T\/}], and
(independently)
P. Arnoux, S. Leli\`evre, and
A. M\'alaga [{\bf ALM\/}],
show that every flat torus, except those
made by identifying the sides of a rectangle,
is realized by some diplotorus.
(Both papers have a separate construction
which works for the rectangular flat tori.)
The 2024 preprint of
F. Lazarus and F. Tallerie [{\bf LT\/}] gives a
universal combinatorial type of triangulation
which does the job simultaneously for
all isometry types. Their triangulation
$1217$ vertices (and $2434$ faces).

There are diplotori, based on a pair of regular
pentagons, which have $10$ vertices.
In 2025,  Vincent Tugay\'e [{\bf Tu\/}],
discovered a $9$-vertex paper torus.
At the same time, one
needs at least $7$ vertices to make a
paper torus because there is
no $6$ vertex triangulation of a torus.

The $7$-vertex triangulation of a
torus, the {\it Moebius triangulation\/},
is unique up to  combinatorial
isomorphism.
The Moebius triangulation gives the famous
embedding of the complete graph $K_7$ in a torus.
Figure 1.1 shows part of the universal cover of the
Moebius torus, as well as a fundamental domain for it.

\begin{center}
\resizebox{!}{2.3in}{\includegraphics{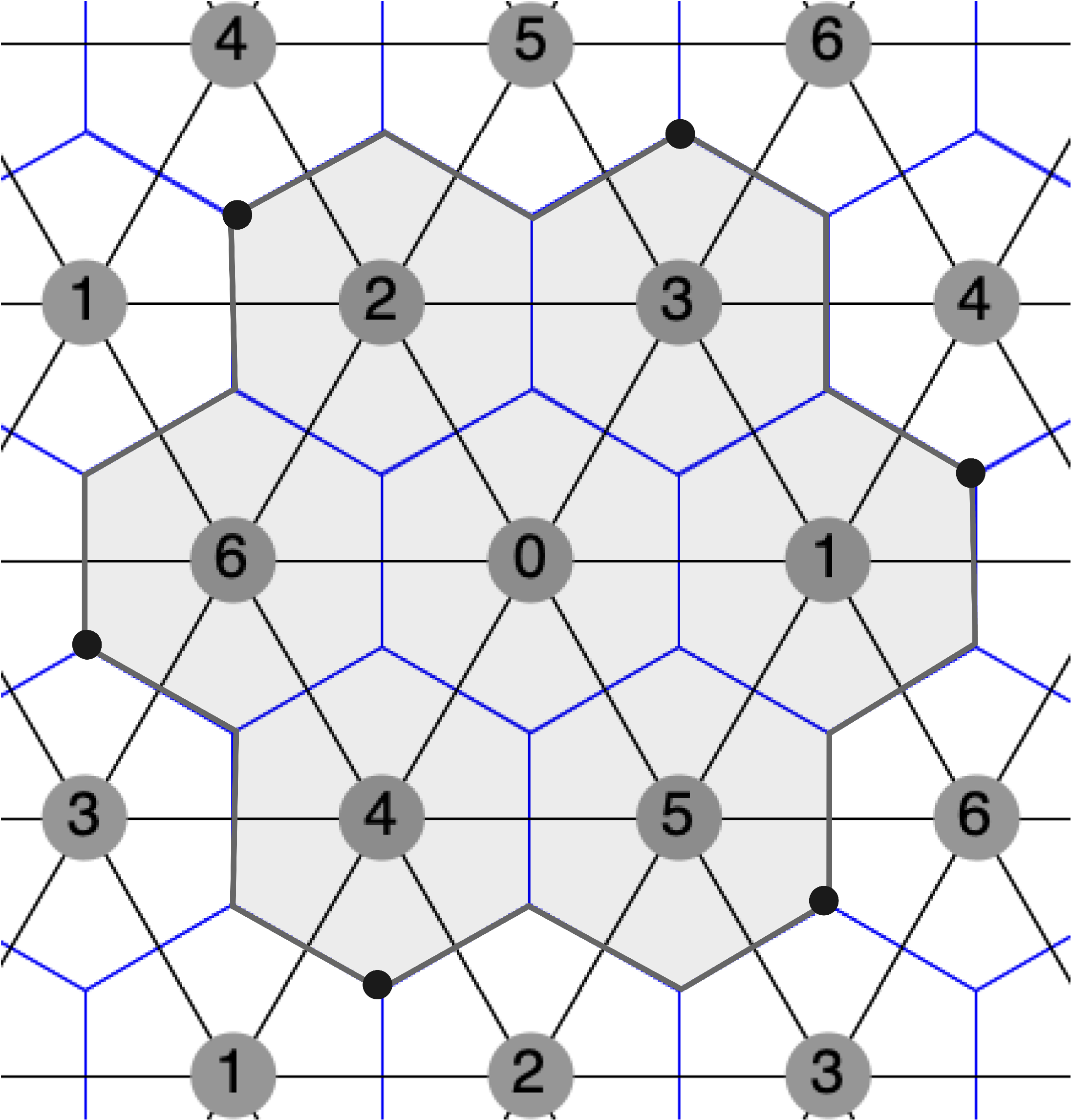}}
\newline
{\bf Figure 1.1:\/} Part of the universal cover of the Moebius triangulation
\end{center}

In $1949$,  Á.  Cs\'asz\'ar  [{\bf Cs\/}] showed
that the Moebius torus has a polyhedral embedding.
However, this embedding is not flat.
See e.g. [{\bf G\/}],  [{\bf HLZ\/}], and [{\bf L\/}]
for more detail and a discussion of related topics.
The 2019 Ph.D. thesis of P. Quintanar Cort\'es,
[{\bf QC\/}] makes significant progress towards
showing that there is no $7$-vertex paper torus.

The results above leave open the question as to the minimum
number of vertices needed to make a paper torus.
The results in this paper settle the question.

\begin{theorem}
  \label{one}
  There does not exist a $7$-vertex paper torus.
\end{theorem}

\begin{theorem}
  \label{two}
  There exists an $8$-vertex paper torus.
\end{theorem}

\subsection{Seven Vertices}

The proof of Theorem \ref{one} is actually rather easy.
Here is the key step.

\begin{lemma}[Hull]
  Suppose $\Omega$ is an embedded $7$-vertex
  polyhedral torus, not necessarily flat.
  Let $H$ be the convex hull of $\Omega$.
  Suppose $\Omega$ has all $7$ vertices in $\partial H$.
  Then $\Omega$ has at least one
  vertex $P$ such that all $6$ triangles of $\Omega$ incident to $P$
  are in $\partial H$.
\end{lemma}

I will deduce from the Hull Lemma that
a paper torus must have a vertex $P$ contained
in the interior of the convex hull.  I will then use
Crofton's formula to show that the cone angle
at $P$ exceeds $2\pi$.  This contradiction finishes
the proof of Theorem \ref{one}.

The Hull Lemma is a consequence of a
stronger result contained in the
1991 paper [{\bf BE\/}] of
J. Bokowski and A. Eggert.
Bokowski and Eggert use oriented matroids to classify
the combinatorial types of $7$-vertex polyhedral
tori. One consequence of their classification is the
{\it Hull Theorem\/}:
there are no embedded $7$-vertex polyhedral tori at all
having all $7$ vertices on their convex hull boundary.
(See the remarks after [{\bf BE\/}, Theorem 3.7]. This
is my name for their result.)
In my opinion, had these authors asked the
flatness question for $7$-vertex polyhedral tori, they
most likely would have been
able to answer it.

In order to keep this paper self-contained,
I will give a light and entirely combinatorial
computer-assisted proof of the
Hull Lemma.  It is not very hard to promote the
Hull Lemma to the Hull Theorem, but since this
is not needed for the proof of Theorem \ref{one}
I will not do it.  I give the argument in my
supplementary notes
[{\bf S\/}].

\subsection{Eight Vertices}

Motivated by our proof of Theorem \ref{one}, I 
tried to rule out the existence of $8$-vertex paper
tori having all $8$ vertices on their convex hull
boundary.  This did not work, so I decided to find
these kinds of examples and study them.

There are $7$ combinatorial types of triangulation
of the torus having $8$ vertices.  Of these $7$ there
is (in my opinion) a {\it worst triangulation\/} and
a {\it best triangulation\/}.   The worst triangulation
is obtained from the Moebius triangulation by
subdividing one of the triangles.  The best triangulation
has degree sequence $66666666$ and is vertex-transitive.
This triangulation is very much like the Moebius
triangulation.  Figure 1.2 below shows part of the universal
cover of this triangulation, as well as a fundamental
domain.  I will explain the meaning of the blue
triangles below.  The triangles all have the form
$(a,a+1,a+3)$ and $(a,a+2,a+3)$ with indices
taken mod $8$.

\begin{center}
\resizebox{!}{2.8in}{\includegraphics{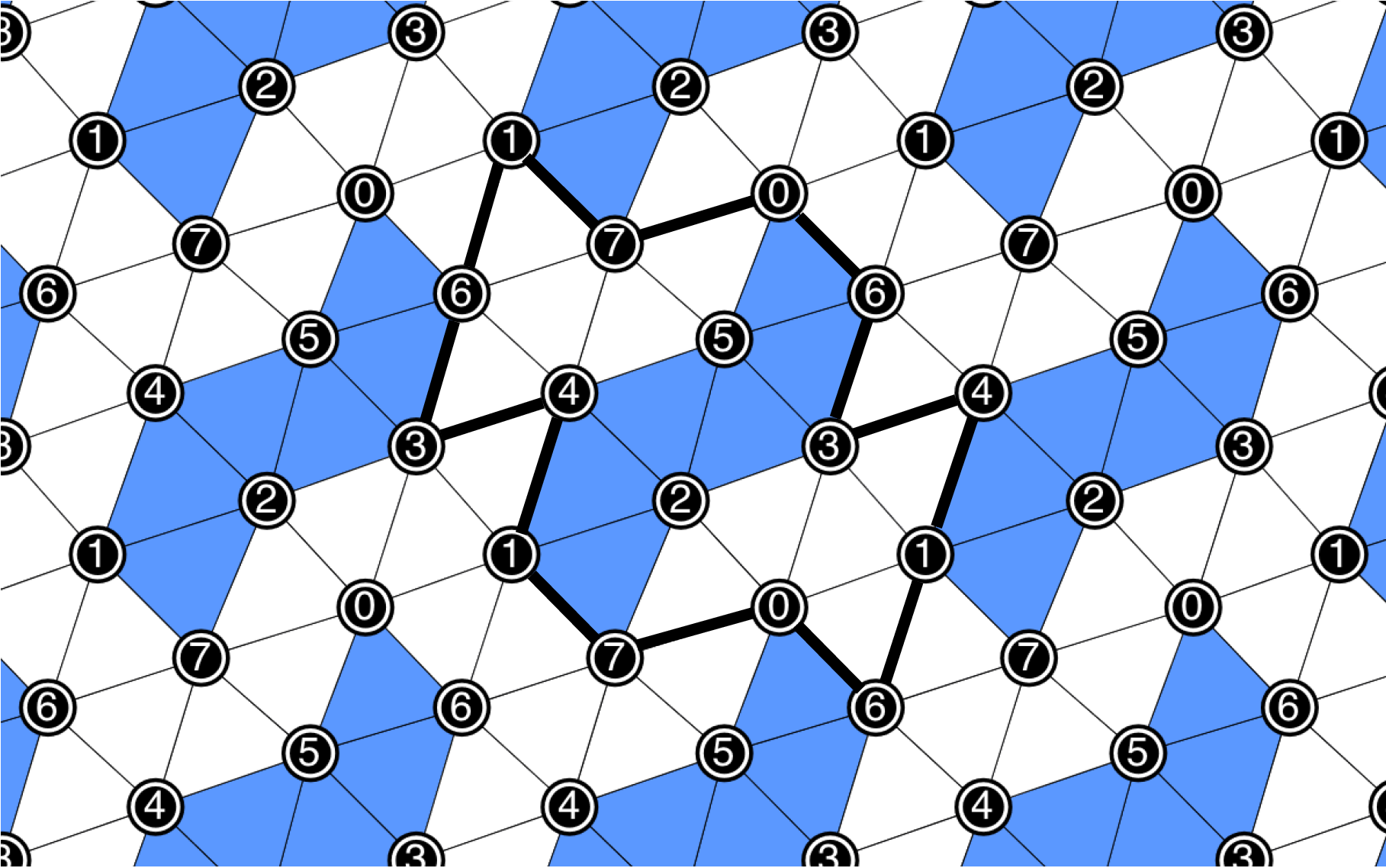}}
\newline
{\bf Figure 1.2:\/} The universal cover of the best $8$ vertex
triangulation
\end{center}

After an intense period of experimentation, I discovered some
$8$-vertex paper tori.
In \S \ref{discuss} I will explain in detail how I found these.
In the example below all the coordinates except
$z_0,z_1,z_2$ are exact.  The
listings of $z_0,z_1,z_2$ are $32$-digit truncations of
the true values.

\begin{equation}
  \label{pup}
\begin{matrix}
+0.64 & -0.20 & 1 \\
-1.09 & +0.38 & z_1 \\
-0.25 & +0.51 & z_2 \\
+0.78 & +0.62 & z_3 \\
-0.78 & -0.62 & z_3 \\
+0.25 & -0.51 & z_2 \\
+1.09 & -0.38 & z_1 \\
-0.64 & +0.20 & 1
\end{matrix}
\hskip 20 pt
\begin{matrix}
\begin{aligned}
z_1 &= 0.0206\ 6632\ 6669\ 8443\ 6159\ 8992\ 3371\ 8861 \\
z_2 &= 0.0048\ 5312\ 7706\ 5192\ 8720\ 4090\ 7479\ 6169 \\
z_3 &= 0.0082\ 2752\ 1455\ 6137\ 1645\ 5791\ 2547\ 8661 
\end{aligned}
\end{matrix}
\end{equation}

Following a suggestion of Peter Doyle, I call
this example and ones like it 
{\it pup tents\/}, on account of their appearance.
See Figure 4.1.   Specifically, a pup tent is an
$8$-vertex paper torus with the following
properties:
\begin{itemize}
\item It is based on the $8$-vertex uniform triangulation.
\item It has order $2$ symmetry given by the map
  $(x,y,z) \to (-x,-y,z)$.
\item It has $6$ faces in its convex hull boundary,
  corresponding to the blue faces in Figure 1.2.
\end{itemize}

Here I discuss the proof that the example in
Equation \ref{pup} is really the truncation of
an $8$-vertex paper torus.
If we leave $z_1,z_2,z_3$ unspecified, the
example in Equation \ref{pup} is part of a $3$-parameter
family of examples we call $T(z_1,z_2,z_3)$.
Let $\theta_k(z_1,z_2,z_3)$ be the cone angle of
$T(z_1,z_2,z_3)$ at vertex $k$.
Let
\begin{equation}
  \label{proof8}
  F(z_1,z_2,z_3)=(\theta_1,\theta_2,\theta_3).
  \end{equation}
If $F(z_1,z_2,z_3)=(2\pi,2\pi,2\pi)$
and $T(z_1,z_2,z_3)$ is embedded then we have our proof.
Let $dF$ be the Jacobian of $F$, computed at our example.
Here we show $dF$ with $4$ digits of accuracy.

  \begin{equation}
    \label{JACOB}
    \left[\begin{matrix}
    -1.6627 \ldots &	1.4699\ldots &	-0.0244 \ldots \\
     1.4699 \ldots & 	0.04714 \ldots &	-1.1300 \ldots \\
     -0.0244 \ldots &	-1.1300 \ldots &	1.9196 \ldots
   \end{matrix}\right]
 \end{equation}
 The matrix $dF$ is an invertible matrix, and all of its eigenvalues exceed
   $3/10$ in absolute value.  This means that the inverse
   $dF^{-1}$ is also small.   In other words, $dF$ is a
   well-conditioned matrix.

   Our proof is an interplay of three things. First, we show that
  our example is robustly embedded in the sense that moderate
  perturbations keep it embedded.  Second, as we say above, our
  example is very close to being flat.
  Third we show, by bounding derivatives,
  that our Jacobian matrix stays fairly constant and well conditioned
  throughout a healthy neighborhood of our example. We then use an
  effective version of the Inverse Function Theorem so show that
  a true paper torus is quite near our example.    We will
  add enough slack to our proof so that it would work
  for examples which are fairly near the one in
  Equation \ref{pup}.

  Given the flexibility of our argument, we get the
  following easy corollary.

\begin{corollary}
  The space $\cal X$ of pup tents contains an open
  $6$-dimensional ball.
  \end{corollary}

  \startproof
  We start with the example based on Equation \ref{pup} and
  we then perturb our points symmetrically  by changing
  the $X$ and $Y$ coordinates but not the $Z$-coordinates.
  As long as the perturbation is small enough, the Inverse
  Function Theorem selects $Z$-coordinates for us which
  give a pup tent.   We have $8$ free parameters to play with
  and then these are grouped into a $2$-dimensional
  foliation under the action of similarities of the $XY$-plane
  which fix the origin.  A cross-section of this foliation gives
  us our $6$-dimensional open manifold of inequivalent pup tents.
  \endproof
  
  \subsection{Further Developments}

  Since I wrote the first version of this paper,
  there have been a number of developments.
  First of all, some of us at an ICERM workshop in
  August 2025 printed out pup tents from a
  paper model I made (with advice from
  Alba M\'alaga).    With some amount of
  dexterity, you can use the paper model
  in [{\bf S\/}] to make your own.  This model
  is based on a pup tent that has slightly
  different coordinates than the one in
  Equation \ref{pup}, but the shape is so
  close you can hardly notice the difference.

Second, Peter Doyle made an amazing discovery of
some degenerate versions of the pup tent which realize
every flat structure.  We call the set of these degenerate
pup tents {\it the golden valley\/}.  Peter and I
then used the golden valley to get some further
resuts about $8$ vertex paper tori.
   Our sequel paper
  [{\bf DS\/}] has many of our findings.
  Here are two example results:
  \begin{enumerate}
    \item Any flat torus which lacks reflection symmetry -- this
      accounts
      for almost all of them -- is realized by a pup tent.
  \item Given any $\epsilon>0$ there is an $8$-vertex paper torus
    within $\epsilon$ of the unit equilateral triangle in the
    Hausdorff metric.
  \end{enumerate}
  Also, the proofs in [{\bf DS\/}] give a second and independent proof of
  Theorem \ref{two}.

  Third, after hearing about the pup tent, my undergraduate student,
  Zhengyu (Byron) Zou tried for similar results in
  hyperbolic geometry.   Recently he found a $10$-vertex
  intrinsically hyperbolic embedding of a genus $2$
  surface in $\H^3$.  This is vertex-minimal, because
  you need $10$ vertices to triangulate a genus $2$ surface.
  See [{\bf Z\/}].
  Even more recently, he discovered a $12$-vertex
  hyperbolic embedding of a genus $2$ surface
  which has order $2$ symmetry.  The order $2$
  symmetry implements the hyperelliptic involution.
  
\subsection{Computer Assistance}

My proof is heavily computer assisted, though
all the important calculations are done with
exact integer arithmetic.
The reader can download all the code for this paper at
{\bf www.math.brown.edu/$\sim$res/Papers/VertexMinimal.tar\/}.
\newline

 I used ChatGPT-4o and ChatGPT-5 in
various ways to help me with this project.
Since this is my first time doing a project with the assistance
of ChatGPT, I think that this warrants some discussion.
I used ChatGPT to discuss high level mathematical ideas,
to find typos and other inconsistencies in my paper,
to format data both for the paper and for my computer code,
and also to document the efficiency of my programs.
There was one astounding thing ChatGPT5 did.
In a stroke of brilliance, it found on its own the
expression for the denominator of $g_{11}$ at the
end of the paper.  I might not have arrived at that.
I did not let
ChatGPT compose prose for me and I ignored its writing advice.
I want to retain $100$ percent of my
own voice when I write something.
Overall, I found ChatGPT to be a helpful and enthusiastic
companion as I worked.

 \subsection{Acknowledgements}
 
  I thank Samuel Leli\`evre and
  Alba M\'alaga for telling me all about
  flat tori (over a period of some years)
  and in particular telling me about the flatness problem
  and supplying me with some historical context.  The
  main thing that inspired me to work on this project
  was the infectious enthusiasm of Samuel and Alba.
  I also thank Peter Doyle,  Fabian Lander,
  Noah Montgomery, Stepan Paul,
 Saul Schleimer, and Byron Zou  for helpful
  discussions about this material. I thank
  Frank Lutz for making his manifolds webpage, which
  lets you download the combinatorial data for
  triangulations of surfaces and other manifolds.
  This page was very useful to me.
  
   I thank the IHES, the
  Hamilton Institute, and the
  Isaac Newton Institute, where I started working
  on this paper.  This work is also
  supported by a Simons Sabbatical Fellowship
  and the Mazur Chair at IHES.

\newpage

\section{Proof of Theorem \ref{one}}

\subsection{Crofton's Formula}

Crofton's Formula is a classic result from
integral geometry. See [{\bf S\/}]. It applies to any rectifiable
arc on the unit sphere $S^2$.     We just
need this result for finite unions of
arcs of great circles, which we call
{\it spherical polygons\/}.  

The space of oriented great circles in $S^2$
is canonically bijective with $S^2$ itself,
and inherits a canonical probability measure $\sigma$.
Given a spherical polygon $\gamma$ and a
great circle $C$, we let
$\#(C \cap \gamma)$ denote the
number of intersection points. (We can ignore the
finitely many great circles for which this is infinite.)
Crofton's formula is as follows.
\begin{equation}
  \label{crofton}
  {\rm length\/}(\gamma) = \pi \int_{\cal S} \#(C \cap \gamma)\ d\sigma.
\end{equation}

At least for spherical polygons, this has a swift proof.
Both sides of Equation \ref{crofton} are additive, so it
suffices to prove the equation for a spherical polygon
that is just an arc of a single great circle.  By continuity,
it suffices to prove the result for great circle arcs whose
length is a rational multiple of $2\pi$.  By additivity once
again,  it suffices to prove Crofton's formula
for the great circle itself.  But then the formula is obvious.

Here is a well-known consequence of Crofton's formula.

\begin{lemma}
  \label{hemi}
Suppose $\gamma$ is a spherical polygon which is also a topological loop.
If
$\gamma$ has length at most $2\pi$ then
$\gamma$ lies in some hemisphere of $S^2$.
\end{lemma}

\startproof
Crofton's formula says that some great circle $C$
intersects $\gamma$ at most once.
But then $\gamma$ lies in one of the two
hemispheres defined by $C$.
When $\gamma$ has length exactly $2\pi$
we can shorten $\gamma$ a bit by cutting
a small corner off.  The shortened curve
then lies in a hemisphere.  Taking a limit
as the cut corner tends to $0$ in length,
we see that $\gamma$ itself lies in a
hemisphere.
\endproof

\subsection{Proof of the Main Theorem}

Now suppose that $\Omega$ is a flat embedded
$7$-vertex polyhedral torus.  Let $H$ be the
convex hull of $\Omega$.
Call a vertex of $\Omega$ {\it interior\/} if it
does not lie in $\partial H$.  Given any
vertex $P$ of $\Omega$ let $\Omega_P$ denote the
union of $6$ triangles of $\Omega$ incident to $P$.

\begin{lemma}
  $\Omega$ has an interior vertex.
\end{lemma}

\startproof
Suppose that $\Omega$ does not have an interior vertex.
By the Hull Lemma, $\Omega$ has a vertex $P$ such that
all $6$ triangles of
$\Omega$ incident to $P$ lie in
$\partial H$.  Let $\Omega_P$ be the union of triangles
incident to $P$ as above.

A {\it support plane\/} for $H$ is a plane which intersects
the boundary $\partial H$ and is disjoint from the interior of $H$.
Since the cone angle of $\Omega$ at $P$ is $2\pi$, and since
$\Omega_P \subset \partial H$, we see that $H$ cannot be strictly
convex at $P$.  Hence  some support plane
through $P$ intersects $\Omega_P$ in (at least) a line segment
through $P$ that contains $P$ in its relative interior.
But then the intersection of $\partial H$ with a small
sphere centered at $P$ is a spherical polygon consisting
of two arcs which connect the same pair of antipodal
points.

Since the total length of $\gamma$ is $2\pi$, and
since a great semicircle is the unique path on
$S^2$ of length at most $\pi$ connecting two antipodal points,
we see that $\gamma$ must be a union of two great semicircles.
Hence $\Omega_P$ is contained in the union of two planes
whose intersection contains $P$.

One of these two planes, $\Pi$, contains at
least $3$ consecutive triangles of $\Omega_P$.
But then $\Pi$ contains at least $5$ vertices of
$\Omega$.  Since the $1$-skeleton of the triangulation of
$\Omega$ is the complete graph $K_7$, we see that $\Pi$ contains
an embedded copy of $K_5$.
This contradicts the fact that $K_5$ is not planar.
\endproof

Now we know that our flat embedded $7$-vertex torus
$\Omega$ has an interior vertex $P$.

\begin{lemma}
  \label{cone}
  The cone angle $\theta$ at $P$ exceeds $2\pi$.
\end{lemma}

\startproof
We translate so that $P$ is the origin.  Let $S^2$ be
the unit sphere. Let $\theta_P$ denote the cone angle
of $\Omega$ at $P$.  Let $\ell(\cdot)$ stand for length.

Let $\Omega_P$ be the union of $6$ triangles of
$\Omega$ incident to $P$, as above.
Let $\widehat \Omega_P$ denote the union of rays emanating
from $P$ whose initial portions are contained in
$\Omega_P$.    Let $L_P=\widehat \Omega_P \cap S^2$.
We have $\theta_P=\ell(L_P)$.

Since $P$ is in the interior of the convex hull of
$\Omega$ and since $\widehat \Omega_P$ contains all the vertices
of this convex hull, we see that $\widehat \Omega_P$ cannot lie in
any halfspace bounded by a plane through the origin.
This means that $L_P$ cannot lie in any
hemisphere of $S^2$.  By Lemma \ref{hemi},
we have $\ell(L_P)>2\pi$. Hence $\theta_P>2\pi$.
\endproof

This completes the proof of the Main Theorem
modulo the Hull Lemma.  We mention again
that the Hull Lemma is an immediate consequence
of (what I call) the Hull Theorem in
[{\bf BE\/}].    Readers who are keen to
see the existence of the $8$-vertex paper
torus can stop reading now.  Readers
who want a self-contained proof of the
Hull Lemma should keep reading.

\subsection{Proof of the Hull Lemma}
\label{combinatorics}

We suppose that
$\Omega$ is an embedded $7$-vertex polyhedral torus having
$7$ points on its convex hull $H$.   We
can perturb so that the points are in
 general position.  This makes $H$ into a triangulated
 solid polyhedron with $7$ vertices, $15$ edges,
 and $10$ faces. One property we will use repeatedly 
  is that $\Omega$ is
  {\it neighborly\/}:  Every two vertices of $\Omega$ are
  joined by an
 edge in the triangulation.  In particular, every edge of
 $H$ is an edge of the triangulation of $\Omega$.

We call the $15$ edges in $\partial H$ the
{\it external edges\/}.  We call the remaining $6$
edges {\it internal edges\/}.  We say that an
{\it internal edge pattern\/} is a choice
of $6$ distinguished edges from the
$1$-skeleton of the triangulation,
normalized (by symmetry) so that the first
internal edge is $(01)$.
There are $\binom{20}{5}=15504$ different
internal edge patterns.

We say that an {\it external triangle\/} of
$\Omega$ is one that lies in $\partial H$,
and an {\it internal triangle\/} is one that does not.
Each internal edge is incident to two
internal triangles.

\begin{lemma}
  \label{nosolid}
  An internal triangle cannot be bounded by
  $3$ external edges.
\end{lemma}

\startproof
Suppose that such a triangle exists. Call it $\tau$. If
all three edges of $\tau$ lie in $\partial H$
and $\tau \not \subset \partial H$ then
$\tau$ separates $H$ into two components,
both of which contain vertices of $\Omega$
in their interior.  Any path in $H$ connecting
two such vertices must intersect $\tau$.
On the other hand $\Omega-\tau$ is path
connected.  This is a contradiction.
\endproof

By Lemma \ref{nosolid},
the internal edge pattern determines the
set of internal triangles and the set of external triangles.
It is worth pointing out a triangle of $\partial H$ is not
necessarily a triangle of $\Omega$.

Since $\Omega$ is neighborly, the degree of
a vertex in the triangulation of $\partial H$ is
just the number of external edges incident
to that vertex, and this is the same as
$6$ minus the number of internal
edges incident to the vertex.  In particular,
since $H$ is a convex polyhedron and
$\partial H$ is triangulated, the 
degree of each vertex in the triangulation
of $\partial H$ is at least $3$.  In other
words, at most $3$ internal edges
are incident to the same vertex of $\Omega$.
We eliminate all internal edge patterns having
fewer than $3$ internal edges incident to some vertex.

Now we consider the remaining internal edge patterns.
For each vertex $(q)$ we have the list
$\{(q,v_i)\ |\ i=1,...,K_q\}$
of $K_q$ external edges incident to $q$.
We write $K=K_q$ and order these $K=K_q$
vertices cyclically according
to the link of $(q)$ in $\Omega$.
The link of $L_q$ of $(q)$ in $\partial H$ is some
permutation $(w_1,...,w_K)$ of $(v_1,...,v_{K})$.

The following lemma lets us deduce the cycle
structure on $H$ from the pattern of
internal edges.

\begin{lemma}[Cycle Rule]
  In order to be a viable candidate for the link
  of $(q)$ in $\partial H$, the link
  $(w_1,...,w_K)$ must satisfy
  two properties.
\begin{enumerate}
\item $(w_i,w_{i+1})$ must be an external edge for all $i$.
  Indices are taken cyclically 
\item The cycle $(w_1,...,w_{K})$ must be a
  dihedral permutation of $(v_1,...,v_K)$.
  \end{enumerate}
\end{lemma}

  \startproof
The necessity of Condition 1 is obvious.
Condition 2 requires some explanation.
Let $P$ be vertex $(q)$.
Consider the picture in $\partial H$ at $P$.
  Since the points are in general position, $\partial H$
  is a proper convex cone near $P$.  We let
  $\Pi$ be a plane parallel to a support plane through
  $(q)$ that just cuts off a small corner of $H$ near $(q)$.
  Consider the intersection $\partial H \cap \Pi$.
  This is a convex $K$-gon $\Delta$, and the cyclic order
  of $\Delta$ is given by
  $(w_1,...,w_{K})$.

\begin{center}
\resizebox{!}{1.8in}{\includegraphics{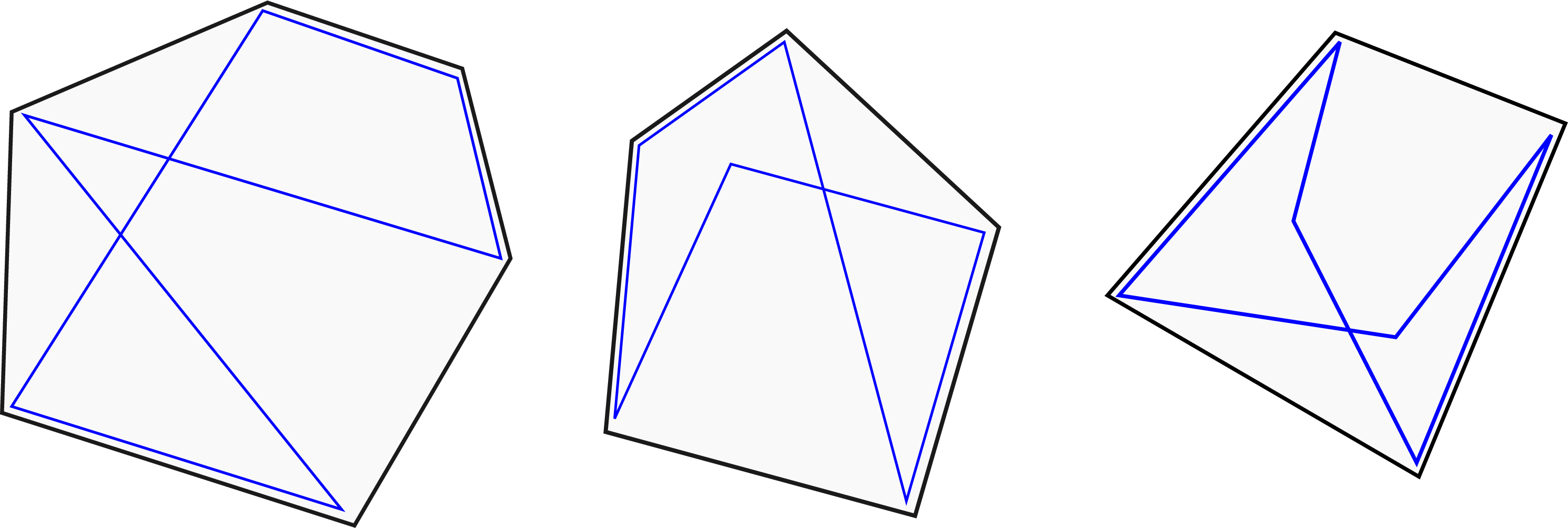}}
\newline
{\bf Figure 2.1:\/}  $\Delta$ (black) and $\gamma$ (blue)
\end{center}

At the same time, and using the notation from
  the proof of Lemma \ref{cone}, let
  $\gamma=\widehat \Omega_P \cap \Pi$.
 The polygonal loop $\gamma$ is contained in the region bounded
  by $\Delta$ and visits the vertices of $\Delta$
  in the order $(v_1,...,v_K)$.  Figure 2.1 shows some
  examples, with $\gamma$ slightly moved off
  $\Delta$ to make the drawing more clear.
  Note that if $K<6$ then there
  will be some extra vertices to $\gamma$ as well;
  this does not matter in the argument.
  If the permutation is
  not dihedral then $\gamma$ cannot
  be an embedded loop.  But then $\Omega_P$ is not embedded,
  a contradiction.
  \endproof

Using the Cycle Rule  (and our computer code)
we eliminate all remaining internal edge patterns except $6$.
These $6$ are the same up to combinatorial isomorphism.
Figure 2.2 shows one of them. 
Evidently, all $6$ triangles incident to vertex $(2)$
are external triangles.  This completes the proof of
the Hull Lemma.

\begin{center}
\resizebox{!}{4in}{\includegraphics{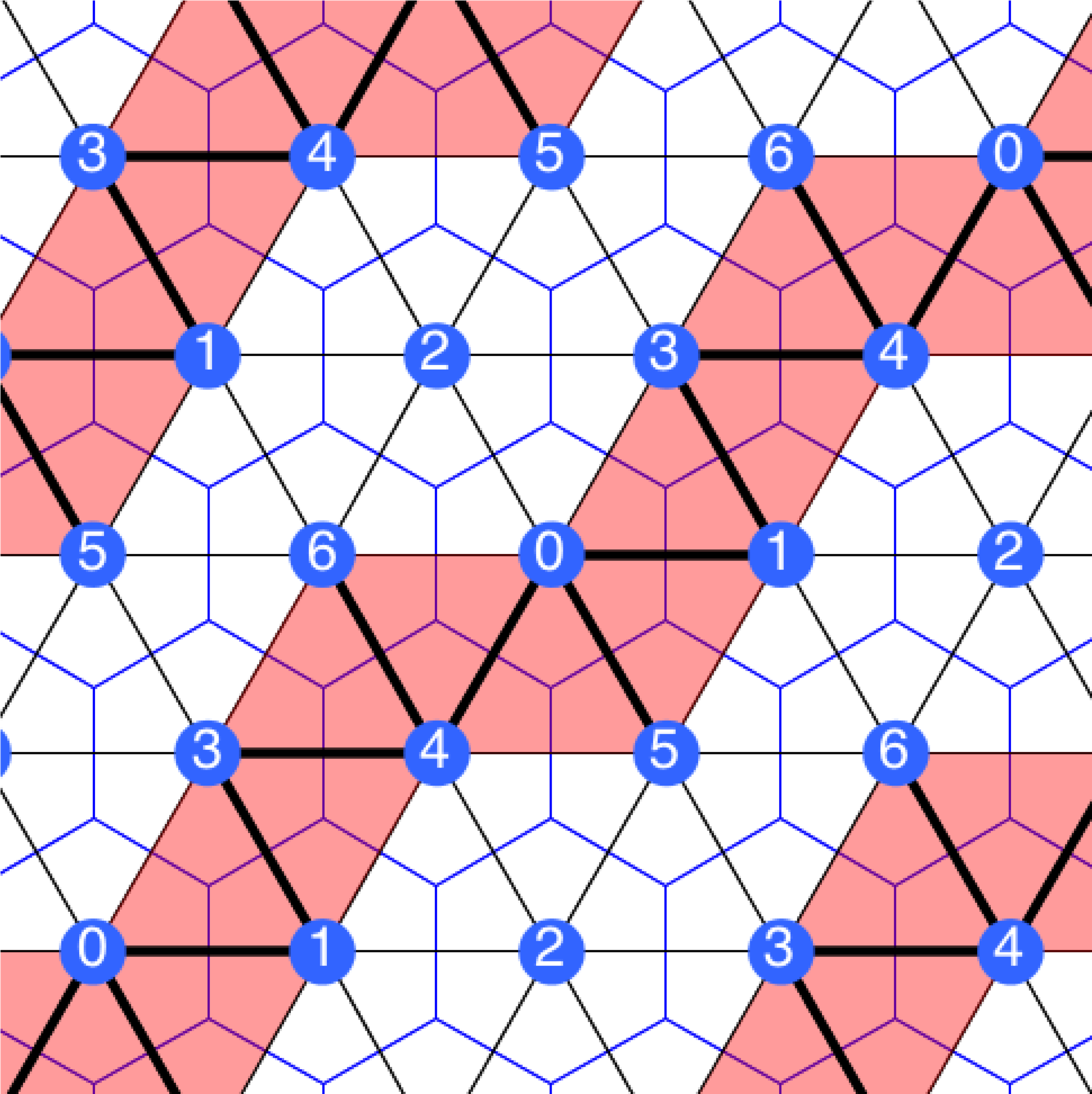}}
\newline
{\bf Figure 2.2:\/}  The one remaining pattern
\end{center}

\noindent
{\bf Remark:\/}
As mentioned in the introduction, we have a separate
geometric argument that rules out the $6$ patterns
like the one in Figure 2.2, thereby giving a second
proof of the Hull Theorem.   See my supplementary
notes [{\bf S\/}] for this argument.

\newpage

\section{Proof of Theorem \ref{two}}

Throughout this chapter $\|\cdot\|$ denotes
the Euclidean norm and $\| \cdot \|_{\infty}$ denotes
the maximum of the absolute values of vector or matrix entries.
Also, $\angle$ denotes angles.   

\subsection{Robust Embedding}
\label{robustproof}

Let $T=T(z_1,z_2,z_3)$ be the example from Equation \ref{pup}.
Let $T(z'_1,z'_2,z'_3)$ be as in Equation \ref{pup}, but
with variables $z_1',z_2',z_3'$.
Let
\begin{equation}
  |T'-T|= \max(|z_1'-z_1|,|z_2'-z_2|,|z_3'-z_3|).
\end{equation}

We say that $T$ is $K$-{\it robustly embedded\/} if and only if
$T'$ is embedded whenever
$|T'-T|\leq K$.  In this section we prove the following result.

\begin{lemma}[Robust Embedding]
  $T$ is $10^{-4}$ robustly embedded.
\end{lemma}

We have $120 = 24 + 72 + 24$ triangle pairs
to consider.  The first $24$ are pairs having
$0$ vertices in common.  The next $72$ are
pairs having $1$ vertex in common. The
last $24$ are pairs having $2$ vertices in common.
We don't need to check these last $24$.

\begin{lemma}
  \label{eee}
  Let $T'$ be a perturbation of $T$ such that
   all triangle pairs with no vertices in
  common are disjoint and all pairs with
  one vertex in common only intersect at
  the common vertex.  Then $T'$ is  embedded.
\end{lemma}

\startproof
We will argue by contradiction.  If
$T'$ is not embedded,
then we can find a pair of triangles
$\tau_1,\tau_2$ having two
vertices $v,w$ in common, which
overlap on an open set.
Let $x_j$ be the third vertex
of $\tau_j$.
After relabeling if necessary,
we can assume that some
initial segment $s$ of the ray
$\overrightarrow{vx_1}$ lies
in $\tau_2$.  But then let
$\beta$ be the triangle of $T'$
that shares the edge $vx_1$ with
$\tau_1$.   The triangle
$\beta$ only has the vertex $v$
in common with $\tau_2$ but
nevertheless intersects $\tau_2$
in a set that contains the
segment $s$.
This is a contradiction.
\endproof

\noindent
{\bf The Disjoint Case:\/}
Suppose $(\Delta_0,\Delta_1)$ is a pair of
triangles with no common vertices.
Given a vector $L$ we define
\begin{equation}
  m_j(L)=\min_{v \in \Delta_j} v \cdot L, \hskip 30 pt
  M_j(L)=\max_{v \in \Delta_j} v \cdot L.
\end{equation}
By convexity, we can compute these quantities just
by computing the dot products on the vertices.
We call $(\Delta_0,\Delta_1)$ $\lambda$-{\it separated\/} if
there exists a vector $L$ with
$\|L\|_{\infty} \leq 1$ such that one of
  the two equations holds:
  \begin{equation}
    \label{sep1}
    M_0(L)+\lambda<m_1(L), \hskip 30 pt
    M_1(L) + \lambda<m_0(L),
  \end{equation}
  We have $\Delta_0 \cap \Delta_1 = \emptyset$ as
  long as this pair is $0$-separated.

  We show by direct calculation that all $24$ pairs
  are $2 \times 10^{-4}$ separated.
  Let $Q$ be the set of integer vectors
  $(x,y,z)$ with $\max(|x|,|y|,|z|)=300$.
  We simply enumerate all vectors in
  the scaled set $\frac{1}{300}Q$
  and find that at least one  works
  for each pair.   We make this an integer
  calculation by scaling up these vectors
  by $300$ and scaling up $T$ by $10^{32}$.
  Our separation constant for this scaled problem is
  $6 \times 10^{30}$.

Let $T'$ be another example with
$|T' - T| \leq 10^{-4}$.  Consider the triangle
pair $(\Delta_0',\Delta_1')$ of $T'$ that corresponds to $(\Delta_0,\Delta_1)$.
In moving from $\Delta_j$ to $\Delta_j'$, no vertex moves
by more than $10^{-4}$, and at most one of the coordinates changes per vertex.
Hence, each dot product above changes by at most $10^{-4}$.
By the triangle inequality,
$(\Delta_0',\Delta_1')$ is $0$-separated.
Hence $\Delta_0' \cap \Delta_1'=\emptyset$.
\newline
\newline
{\bf The Common Vertex Case:\/}
Now we consider when $(\Delta_0,\Delta_1)$ have a single common
vertex $v$.
Let $e_j$ be the edge of $\Delta_j$ opposite $v$.
This time define
\begin{equation}
  m_j(L)=\min_{w \in e_j} w \cdot L, \hskip 30 pt
   M_j(L)=\max_{w \in e_j} w \cdot L.
\end{equation}
This time
we call $(\Delta_0,\Delta_1)$ $\lambda$-{\it separated\/} if
there exists a vector $L$ with
$\|L\|_{\infty} \leq 1$ such that one of
  the two equations holds:
  \begin{equation}
    \label{sep2}
    M_0(L)+\lambda<v \cdot L<m_1(L)-\lambda, \hskip 30 pt
    M_1(L) + \lambda<v \cdot L<m_0(L)-\lambda,
  \end{equation}
  If $(\Delta_0,\Delta_1)$ is $0$-separated then $\Delta_0 \cap
  \Delta_1=\{v\}$.
  We show by exactly the same kind of calculation that
  all $72$ pairs are $2 \times 10^{-4}$-separated.

  But now the same argument as above shows that for
  $T'$ (as above)  all such pairs are $0$-separated.
  Hence $\Delta_0' \cap \Delta_1'=\{v'\}$.
  \newline

  The two cases above combine with Lemma
  \ref{eee} to prove that  $T'$ is embedded when
  $|T'-T|<10^{-4}$.  Hence
  $T$ is $10^{-4}$-robustly embedded.
  \newline
  \newline
  \noindent
  {\bf Remark:\/}
  There is another approach we could have taken which involves
  computing the signs of the volumes of all the tetrahedra
  made from vertices of $T$.  See [{\bf DS\/}].

\subsection{The Proof Modulo a Crude Bound}

Let $G: \R^3 \to \R^3$ be a smooth map.
Let $B_r(p)$ denote the ball of radius $r$
about $p$.   We say that $G$ is
$\lambda$-{\it expansive\/} on $B_r(p)$
if, for all unit vectors $V$ and all
$q \in B_r(p)$, we have
$\|dG_q(V)\|>2\lambda$ and
$\angle(dF_p(V),dF_q(V))<\pi/3$.

  \begin{lemma}[Surjectivity]
    \label{SUR}
    Let $q=G(p)$.
    If $G$ is $\lambda$-expansive on $B_r(p)$ then
    $B_{\lambda r}(q) \subset G(B_r(p))$.
\end{lemma}

\startproof
Let $O$ be the origin and let $B=B_1(O)$.
By scaling and translation, it suffices to prove
that $B \subset G(B)$ when $p=q=O$ and $\lambda=2$.
Let $\alpha$ be any line segment connecting
$O$ to a point on $\partial B$.
Let $s \mapsto \alpha(s)$ be the unit speed
parametrization of $\alpha$.  Let
$\beta=G(\alpha)$. The parametrization
$s \mapsto \beta(s)$ has speed at least
$2$.   Since $\angle(\beta'(t),\beta'(0))<\pi/3$, the
 projection of $\beta$ onto
the line through $\beta'(0)$ moves with
speed greater than $1$.  But then
$\|\beta(1)-\beta(0)\|>1$.
This shows that $G(\partial B)$ lies entirely outside $B$.
But now the angle condition implies that
$G: \partial B \to \R^3-B$ is homotopic to the
identity map in $\R^3-B$.   Hence
$G(\partial B)$ is
the generator of the
homology group $H_2(\R^3-q)$
for any $q \in B$.
This easily implies that $B \subset G(B)$.
\endproof

Let $F(z_0,z_2,z_3)=(\theta_0,\theta_2,\theta_3)$ be the
map from \S \ref{proof8}.  Let $p \in \R^3$
be the point that corresponds to our torus $T$
from Equation \ref{pup}.
We prove the following result in the next section.

\begin{lemma}[Crude Bound]
  Throughout $B_{10^{-4}}(p)$ we have
  $$\bigg{|} \frac{\partial^2 \theta_k}{\partial z_i \partial
    z_j}\bigg{|}<10^{9}, \hskip 30 pt
\forall\  i,j,k \in \{1,2,3\}.
  $$
\end{lemma}

Here is the main corollary.

\begin{lemma}[Expansion]
  $F$ is $10^{-2}$-expansive on $B_{10^{-12}}(p)$.
\end{lemma}

\startproof
Let $V$ be a unit vector.
Define
  \begin{equation}
    \label{JACOB2}
M=    \left[\begin{matrix}
    -1.66  &	+1.46 &	-0.02 \\
     +1.46  & 	+0.04  &	-1.13  \\
     -0.02  &	-1.13 &	+1.91
   \end{matrix}\right]
 \end{equation}
$M$ is a symmetric and its eigenvalues all exceed $3/10$ in
absolute value.  Hence $\|M(V)\|>3/10$.
We will settle for the weaker bound:
\begin{equation}
  \label{bigimage}
  \|M(V)\|> \frac{1}{10}.
\end{equation}

Define ``error matrices'' $E_q$ and $E'_q$ so that
\begin{equation}
  dF_q=M+E_q=dF_p+E'_q, \hskip 30 pt
  E_q=E_p+E_q'.
\end{equation}
Comparing Equations \ref{JACOB} and
\ref{JACOB2} we see that
\begin{equation}
  \label{JACOB3}
  \|E_p\|_{\infty}<\frac{1}{100}.
\end{equation}

We can join $p$ to $q$ by an axis-aligned path of
length at most $\sqrt 3 \times 10^{-12}$.  Hence,
by the Crude Bound Lemma and integration:
\begin{equation}
  \label{small}
\|E'_q\|_{\infty}
<\sqrt 3 \times 10^{-12} \times 10^{9}<\frac{1}{450}
\hskip 10 pt \Longrightarrow \hskip 10 pt \|E'_q(V)\|<\frac{1}{150}.
\end{equation}
By Equation \ref{JACOB3}, Equation \ref{small}, and the
triangle inequality,
\begin{equation}
  \label{nearness}
  \|E_q\|_{\infty} \leq \|E_p\|_{\infty} + \|E'_q\|_{\infty}<
  \frac{1}{100}+\frac{1}{450}<\frac{1}{75}
\hskip 10 pt \Longrightarrow \hskip 10 pt \|E_{q}(V)\|<\frac{1}{25}.
\end{equation}
By Equation \ref{bigimage}, Equation \ref{nearness}, and the
triangle inequality,
\begin{equation}
  \label{expansive1}
\|dF_q(V)\| \geq \|M(V)\|-\|E_q(V)\| > \frac{1}{10}-
\frac{1}{25}>2 \times \frac{1}{100}.
\end{equation}
By Equation \ref{small}, Equation \ref{expansive1}, and some elementary geometry,
\begin{equation}
  \label{expansive2}
\frac{\|E'_q(V)\|}{\|dF_p(V)\|}<\frac{1/150}{1/50}=\frac{1}{3}<\frac{\sqrt 3}{2}
\hskip 10 pt
\Longrightarrow  \hskip 10 pt \angle(dF_p(V),dF_q(V))<\frac{\pi}{3}.
\end{equation}
Equations \ref{expansive1} and \ref{expansive2} give the desired
bounds.
\endproof

Let $q=F(p)$ and let
$q^*=(2\pi,2\pi,2\pi)$.  We compute
directly, in Mathematica, that $\|q-q_*\|<10^{-31}$.
  By the Surjectivity and Expansion Lemmas,
  $q^*=F(p^*)$ with $\|p-p^*\|<10^{-29}$.
The torus $T^*$ corresponding to
$p^*$ is flat.  Since $T$ is
$10^{-4}$ robustly embedded, we see that
$T^*$ is embedded. The dihedral
angles of $T$ all lie in $[.005,\pi-.005]$, so the tiny change
we make in going from $T$ to $T^*$ keeps $T^*$ a pup tent.
This completes the proof of
Theorem \ref{two} modulo the Crude Bound Lemma.
\newline
\newline
{\bf Remark:\/} All we used for our proof is that
$T$ is $10^{-12}$ robustly embedded and
that $\|q-q^*\|<10^{-14}$.  One can get these
weaker bounds from ordinary floating point
operations in Java.

\subsection{Proof of the Crude Bound Lemma}

We keep the notation from above.
Note that when we vary
$(z_1,z_2,z_3)$, we also vary
$(z_6,z_5,z_4)$ in a symmetric way.
Thus,
$$\frac{\partial \theta_j}{\partial z_1} = \partial_1 \theta_j +
\partial_6 \theta_j, \hskip 20 pt
\frac{\partial^2 \theta_j}{\partial z_1\partial z_2} =
(\partial_{12}+\partial_{62} + \partial_{15} +
\partial_{65})\theta_j.$$
$\partial_i \theta_j$ denotes the rate of change of
$\theta_j$ when we vary {\it only\/}  $z_i$ and
$\partial_{ij} \theta_k$ denotes the rate of change of
$\partial_j \theta_k$ when we vary only $z_i$.
Here we are leaving the realm of $2$-fold
symmetric configurations and just considering
$8$ points in space.

Each $\theta_j$ is a $6$-term expression of
terms of the form $\theta_{ijk}$, where
$\theta_{ijk}$ is the angle of triangle $(i,j,k)$ at
vertex $i$.
Hence, $\partial_{ij} \theta_j$ is a $24$-term sum of
expressions of the form
$\partial_{ab}\theta_{cde}$.
These expressions vanish when
$\{a,b\} \not \subset \{c,d,e\}$.
We call this expression {\it friendly\/} if
$c \not \in \{a,b\}$, and otherwise {\it unfriendly\/}.

\begin{lemma}
  Every unfriendly term is a sum of at most $2$ friendly terms.
\end{lemma}

\startproof
Consider terms of
$\partial_{aa}\theta_{abc}$.
Using the relation
$\theta_a+\theta_b+\theta_c=\pi$ we have
$$\partial_{aa}\theta_{abc}+\partial_{aa}\theta_{bca}+\partial_{aa}\theta_{cab}=0.$$
The terms on the right are friendly.

Given that $\partial_{ab}=\partial_{ba}$, the
only remaining unfriendly terms we have to treat are
those of the form $\partial_{ba}\theta_{abc}$ and
$\partial_{ba}\theta_{acb}$.
We will treat the first one. The second one has the
same treatment.

We also have the relation
 $(\partial_a+ \partial_b + \partial_c) \theta_{abc}=0$.
Geometrically, the variation we are considering simply
translates the triangle up or down along the $Z$-axis
and changes none of the angles.
Applying $\partial_b$ to this relation, we have
$$\partial_{ba}\theta_{abc}= -\partial_{bb} \theta_{abc} -
\partial_{bc}\theta_{abc}.$$
Both terms on the right are friendly.
\endproof

In summary, each partial derivative in
the Crude Bound Lemma is 
a sum at most $48$ friendly term.
To prove the Crude Bound Lemma it suffices to
show that each friendly term is less than
$10^{7}$ in absolute value throughout the
ball $B_{10^{-4}}(p)$.

We say that a torus is {\it relevant\/} if it
corresponds to a point in our tiny ball $B$.
We say that a pair $(V_1,V_2)$ of vectors is
{\it relevant\/} if it is obtained from a
triangle $(P_0,P_1,P_2)$ of a relevant torus by setting
$V_1=P_1-P_0$ and $V_2=P_2-P_0$.
  Each triangle in a relevant torus
  gives rise to $6$ pairs of relevant vectors.

  Let $\vartheta(V_1,V_2)$ be the angle between
  $V_1$ and $V_2$.
  Each friendly term has the form
\begin{equation}
  \label{expression}
  \frac{\partial^2 \vartheta}{\partial w_i
    \partial w_j}, \hskip 30 pt \vartheta= \arccos\bigg(\frac{V_1 \cdot V_2}{\sqrt{(V_1 \cdot V_1) (V_2 \cdot V_2)}}\bigg)
\end{equation}
The square $g_{ij}$ of the expression on the left side of
Equation \ref{expression} is a rational function.
To prove the Crude Bound Lemma, we just
need to prove for all relevant pairs and all $(i,j)$ that
that $g_{ij}<10^{14}$.

\begin{lemma}
  \label{BBB}
  For
  each relevant pair $(V_1,V_2)$ of vectors,
  $$\|V_1 \times V_2\|^2, \|V_1\|^2,\|V_2\|^2 \in (10^{-1},10).$$
\end{lemma}

\startproof
For $T$ we get the bounds $\|V_k\| \in [0.48,2.4]$
and $\|V_1 \times V_2\| \in [.85,2.4]$.  These various
quantities hardly change at all when we switch to a
different relevant torus and change our coordinates
by less than $10^{-4}$.
   \endproof

   \noindent
   {\bf Case 1:\/}
   Let $W_k=(u_k,v_k,0)$ be the projection of $V_k$ to the $XY$-plane.
   We have
   $$g_{21}=g_{12}=
   \frac{\|W_1 \times W_2\|^4}{\|V_1 \times V_2\|^6} \leq 
   \frac{\|V_1 \times V_2\|^4}{\|V_1 \times V_2\|^6} =
   \frac{1}{\|V_1 \times V_2\|^2}<10.
   $$
   The second inequality comes from Lemma \ref{BBB}.
   \newline
   \newline
   {\bf Case 2:\/}
   With $W_k$ defined as above, we have
   $$g_{11}=\frac{ \bigg(t_1 w_1 w_2 + t_2 w_1 w_2^3 + t_3  w_1^3 w_2 
+ t_4 w_1^2 w_2^2 + t_5 w_1^4 + t_6 w_1^2 + t_7\bigg)^2}{\|V_1\|^8 \|V_1
\times V_2\|^6},$$

$$
\begin{aligned}
t_1 &= 3\,\|W_1\|^2 \,\|W_1\times W_2\|^2, 
&\quad t_2 &= 2\,\|W_1\|^4, \\[6pt]
t_3 &= 3\,(W_1\cdot W_2)^2 \;+\; 3\,\|W_1\|^2 \|W_2\|^2, 
&\quad t_4 &= -6\,\|W_1\|^2 (W_1\cdot W_2), \\[6pt]
t_5 &= -2\,\|W_2\|^2 (W_1\cdot W_2), 
&\quad t_6 &= -\,\|W_1\times W_2\|^2 \,(W_1\cdot W_2), \\[6pt]
t_7 &= \|W_1\|^2 \,\|W_1\times W_2\|^2\,(W_1\cdot W_2). &
\end{aligned}
$$

Write this as $g_{11}=N^2/D$.
Using Lemma \ref{BBB} and the fact that
$|w_1|,|w_2|<1$ we see that
$|N|<10^{7/2}$ and $D>10^{-7}$.
Hence $g_{11}<10^{14}$.
A similar calculation works for
   $g_{22}$ and indeed this case follows from that of $g_{11}$ and
   symmetry.
   \newline

This completes the proof of the
Crude Bound Lemma, and thereby
completes the proof of
Theorem \ref{two}.
An $8$-vertex paper torus exists!

\newpage
  
\section{Discussion and Pictures}
\label{discuss}

\subsection{Experimentation}

Recall that there are $7$ combinatorial types of
triangulation of the torus having $8$ vertices.
Figure 1.2 shows the best of these, the one
with uniform degree.
For all but the worst triangulation, I generated lots of embedded examples
by randomly sampling collections of $8$ points on the unit sphere
and testing the embedding condition.

Next, I used a hill climbing algorithm.  The way
this works is that you want to minimize an objective
function on a space.  You start at a random point and
then evaluate the function.   You then make a random
perturbation and move to the new location if the
function is lower. If you are in a smooth manifold and
the function is smooth and nonsingular, then the random perturbations
often move in the direction of the gradient. In this setting
the hill climbing algorithm somewhat mimics gradient
flow.  This approach works pretty well in low dimensions.

My hill-climbing algorithm started with
random embedded examples on the unit sphere and
repeatedly perturbed them -- not necessarily keeping them
on the unit sphere -- with the goal of minimizing an
objective {\it flatness function\/}.  
The flatness function
computes the maximum deviation of a cone angle from
being $2\pi$.

After a lot of experimentation, I implemented
$4$ tricks that helped immensely.

\begin{enumerate}
\item It was useful to first perturb the points within an ellipsoid
surface that contains them, and then to perturb the
ellipsoid.   This method is guaranteed to keep all vertices
on the convex hull boundary.

\item It was extremely useful to control the
combinatorics on the convex hull boundary.
The convex hull boundary has $12$ faces, some of which --
say $K$ of them -- are also faces in the triangulation of the torus.
Call $K$ the {\it face number\/}.   I noticed that the examples
with large face number never led to flat examples.
So I eventually constrained my search to examples with low face
number.

\item It helped to bound some auxiliary
  quantities, like the minimum angles and dihedral angles, to
  prevent degenerations.  Bounds of about $0.00001$ on the
  dihedral angles worked well.
  
\item It helped to force the examples to have $2$-fold rotational
  symmetry.   This simplification came later. After I found some
  asymmetric examples, Samuel Leli\`evre asked me about symmetric
  examples and I tweaked the program.  The program works more reliably
  after incorporating the symmetry.
\end{enumerate}

\subsection{Pictorial Study}

With all the tricks above in place, and a lot of
playing around,
I found some embedded and $2$-fold symmetric
examples of face number $6$ which are flat up
to $10^{-16}$.   I took the nicest one and simplified the
coordinates, feeding it back into the program to make it
near-flat again.   After some fooling around with the
digits I arrived at the example in Equation \ref{pup}.

Figure 4.1 shows three Mathematica plots of the pup tent,
from various angles.

\begin{center}
  \resizebox{!}{2.08in}{\includegraphics{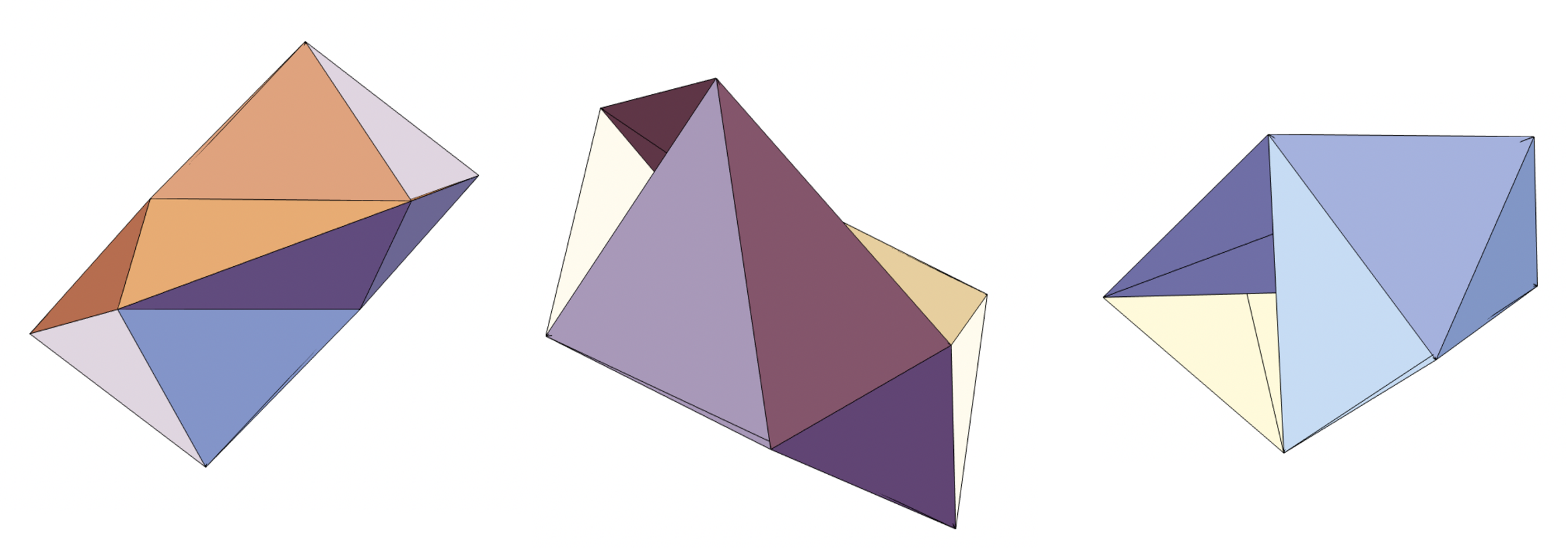}} \newline {\bf
    Figure 4.1:\/} 3D plots of the pup tent
\end{center}

Figure 4.2 shows a numerical computation of the intrinsic geometric
structure of the faces of the pup tent.  We are showing part of the
universal cover of the triangulation.   We have also
highlighted the fundamental domain that corresponds to the
one in Figure 1.2.  The blue triangles correspond to the faces of the
pup tent on the
convex hull boundary. There are $6$ of them in the fundamental domain.

\begin{center}
  \resizebox{!}{2.32in}{\includegraphics{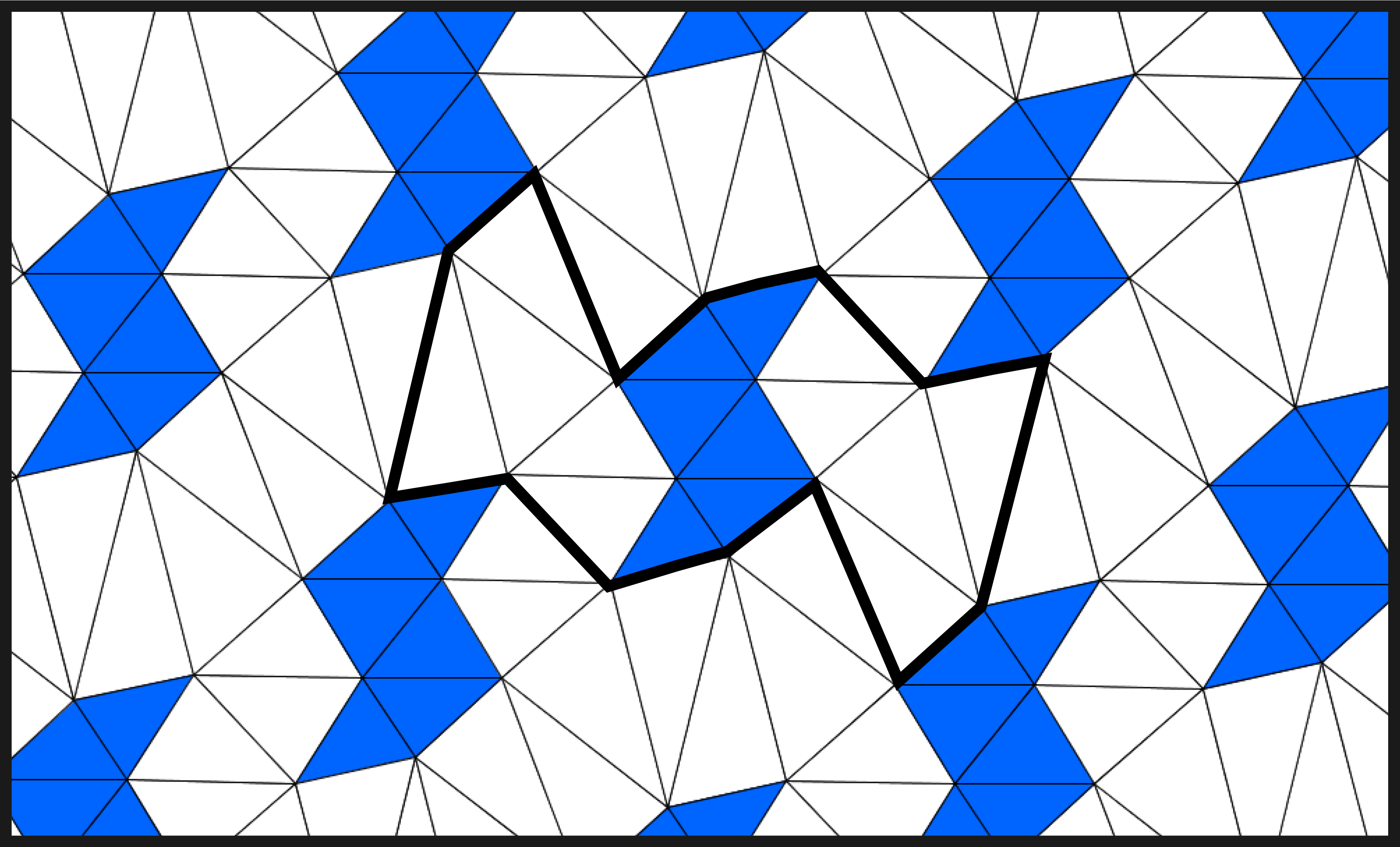}} \newline {\bf
    Figure 4.2:\/} The intrinsic structure.
\end{center}

Figure 4.3 shows a second coloring of the same
fundamental domain.  We
have labeled the points. We list out numerical approximations
to the $8$ blue triangles.  We get the pink ones by
reflecting the blue ones in the origin and applying the
permutation $j \to 7-j$ to the labels.

\begin{center}
  \resizebox{!}{3.6in}{\includegraphics{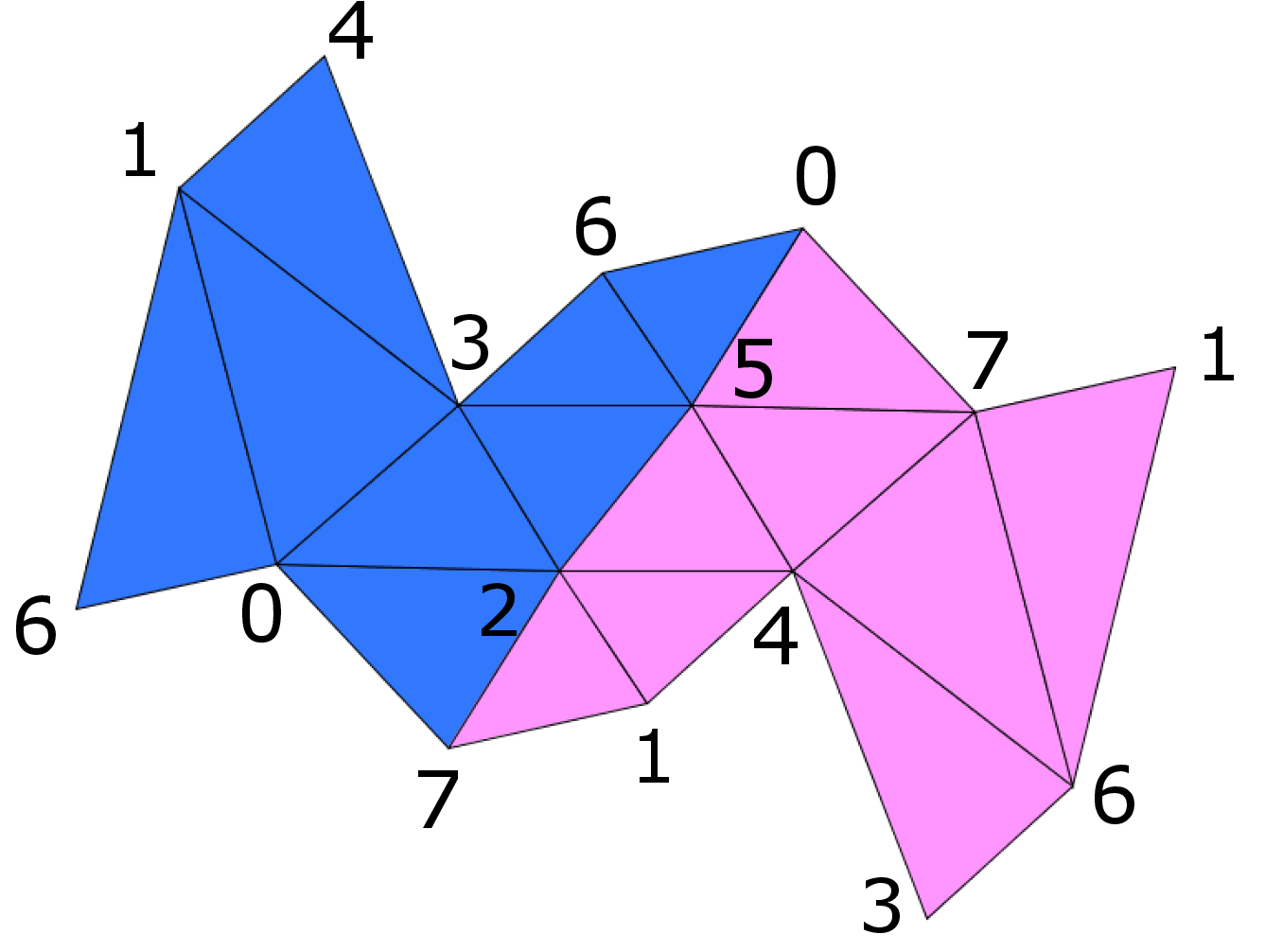}} \newline {\bf
    Figure 4.3:\/} A fundamental domain for the tiling
\end{center}

{\scriptsize
\renewcommand{\arraystretch}{1.3}
\[
\begin{array}{lll}
\{6,0,1\}\!: & \{-2.9340654093 - 0.6455398803 i,\; -1.8673131665 - 0.4088246583 i,\; -2.3837468425 + 1.5965906856 i\}, \\[2pt]
\{3,1,0\}\!: & \{-0.8925492572 + 0.4429076765 i,\; -2.3837468425 + 1.5965906856 i,\; -1.8673131665 - 0.4088246583 i\}, \\[2pt]
\{3,4,1\}\!: & \{-0.8925492572 + 0.4429076765 i,\; -1.6100587412 + 2.3020419213 i,\; -2.3837468425 + 1.5965906856 i\}, \\[2pt]
\{2,0,7\}\!: & \{-0.3555738883 - 0.4429076765 i,\; -1.8673131665 - 0.4088246583 i,\; -0.9478910867 - 1.3850741341 i\}, \\[2pt]
\{3,0,2\}\!: & \{-0.8925492572 + 0.4429076765 i,\; -1.8673131665 - 0.4088246583 i,\; -0.3555738883 - 0.4429076765 i\}, \\[2pt]
\{3,2,5\}\!: & \{-0.8925492572 + 0.4429076765 i,\; -0.3555738883 - 0.4429076765 i,\; +0.3555738883 + 0.4429076765 i\}, \\[2pt]
\{3,5,6\}\!: & \{-0.8925492572 + 0.4429076765 i,\; + 0.3555738883 + 0.4429076765 i,\; -0.1188611560 + 1.1483589122 i\}, \\[2pt]
\{5,0,6\}\!: & \{ +0.3555738883 + 0.4429076765 i \;+1.8673131665 +0.40882465835 i,\;+0.9478910867 +1.3850741341 i\}
\end{array}
\]
\/}

Figures 4.2 and 4.3 look quite tame, but the extrinsic geometry of
the pup tent is
beautiful and very intricate.  Figure 4.4 shows one of the
slices of the pup tent by a plane parallel to the $Z$ axis.
The solid handlebody bounded by the torus is colored orange.
This slice is very spiky and sort of ``barely embedded'',
but the embeddedness is robust compared to the plot precision.

\begin{center}
\resizebox{!}{2.2in}{\includegraphics{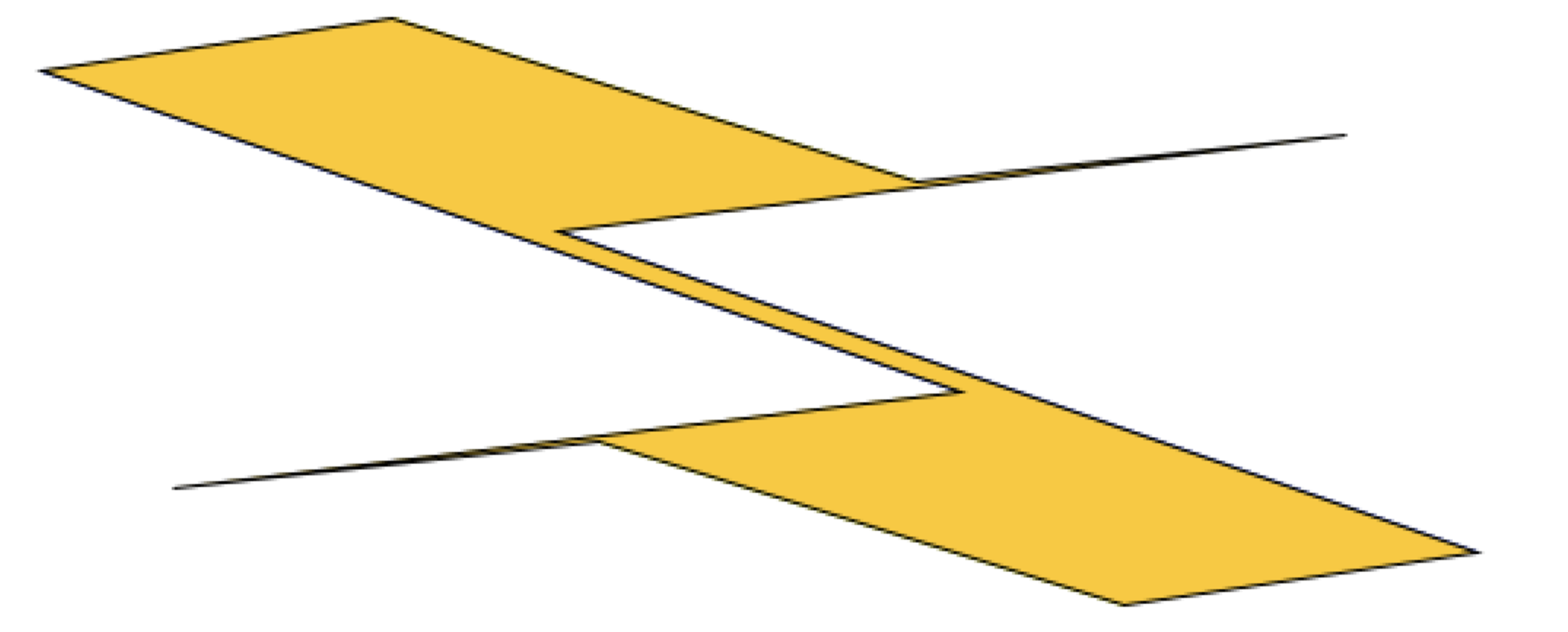}}
\newline
{\bf Figure 4.4:\/} A horizontal slice of the pup tent
\end{center}

Figure 4.5 shows the slice by the $XZ$ plane, which
goes through the axis of symmetry of the torus.
The slice of the handlebody is an annulus.
We have included a closeup of the bottom so that
the reader can see it better.

\begin{center}
\resizebox{!}{4in}{\includegraphics{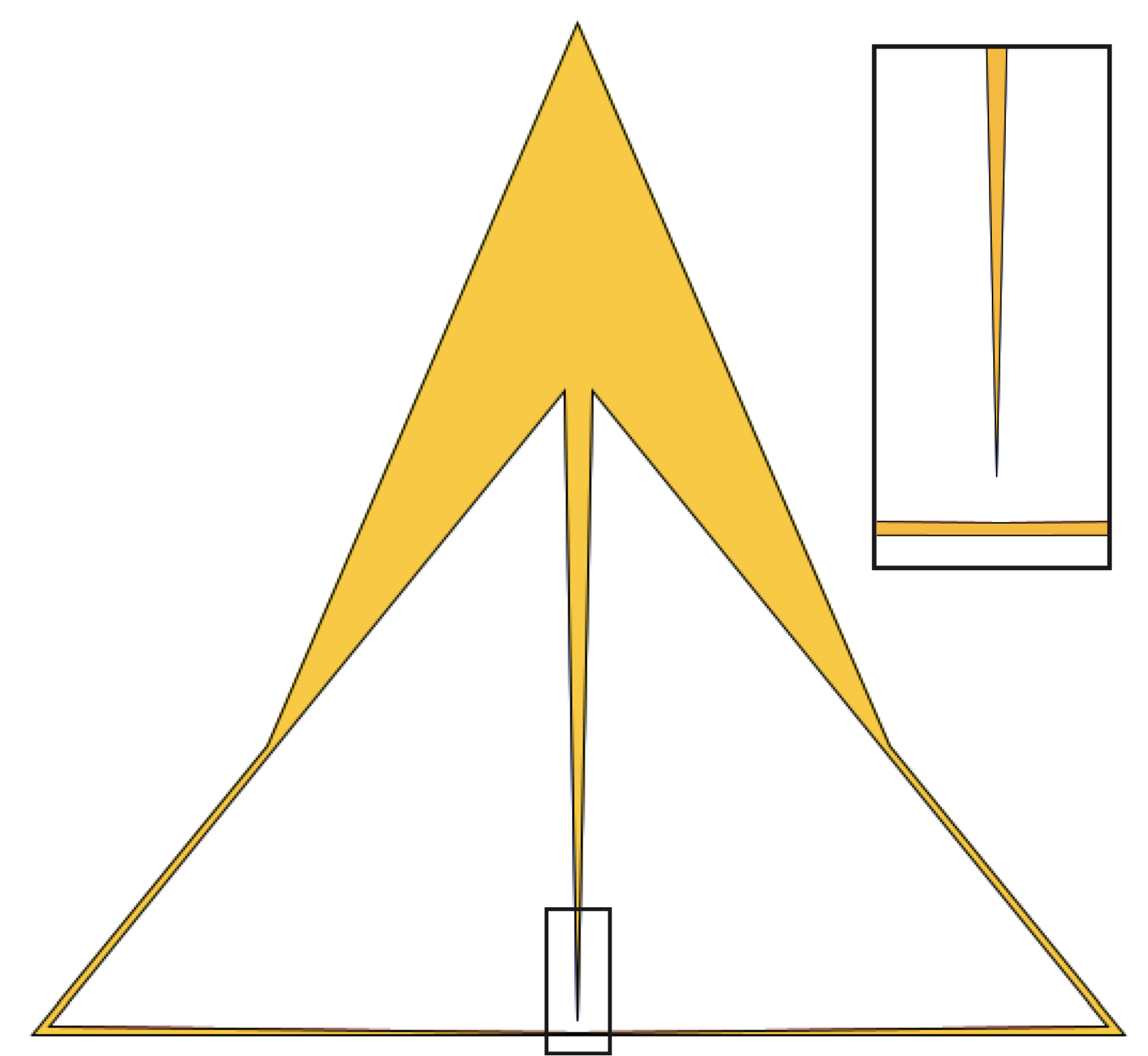}}
\newline
{\bf Figure 4.5:\/}  The $XZ$ planar slice of the pup tent.
\end{center}

Figure 4.6 shows the slice by the $YZ$ plane, with a closeup
of the bottom.

\begin{center}
\resizebox{!}{2.5in}{\includegraphics{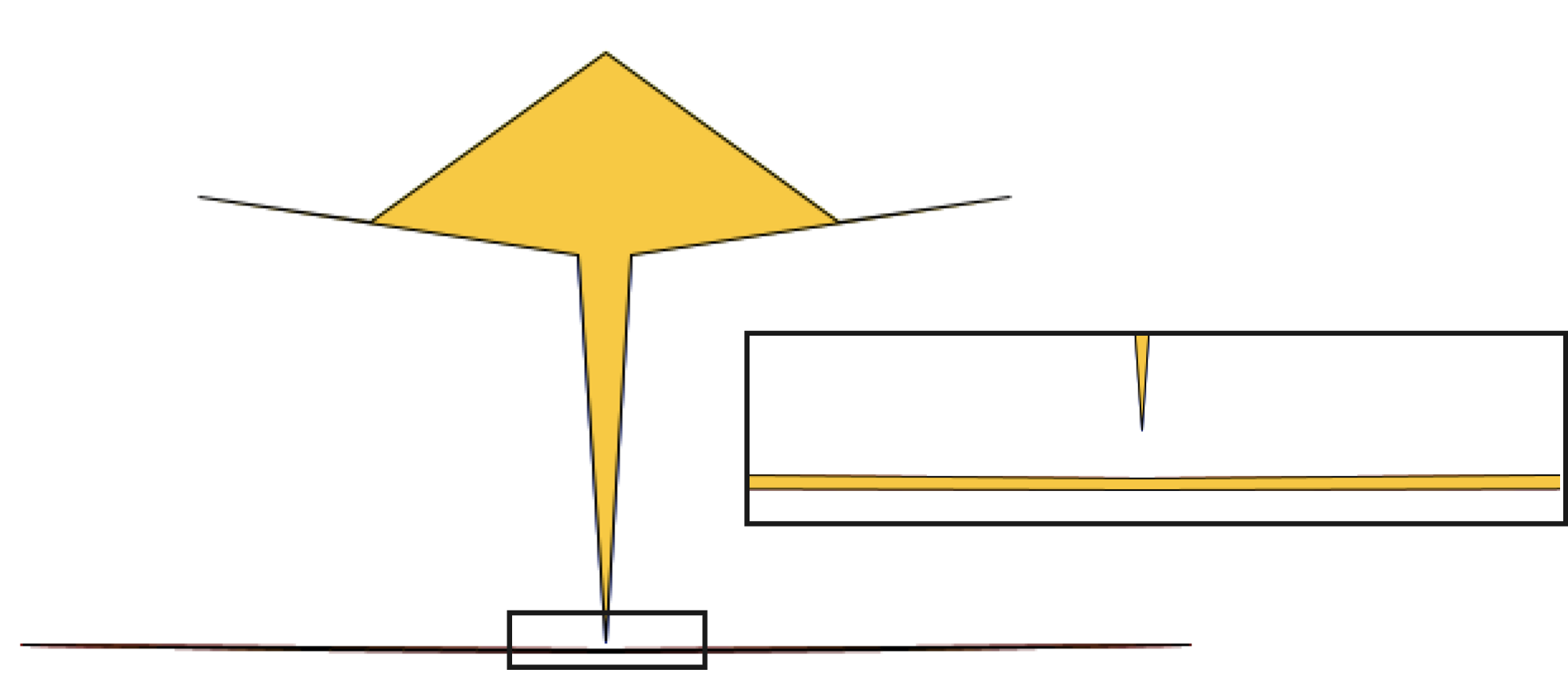}}
\newline
{\bf Figure 4.6:\/} The $YZ$ planar slice of the pup tent
\end{center}

The left side of Figure 4.7 shows a random slice of the pup tent.

\begin{center}
\resizebox{!}{2.4in}{\includegraphics{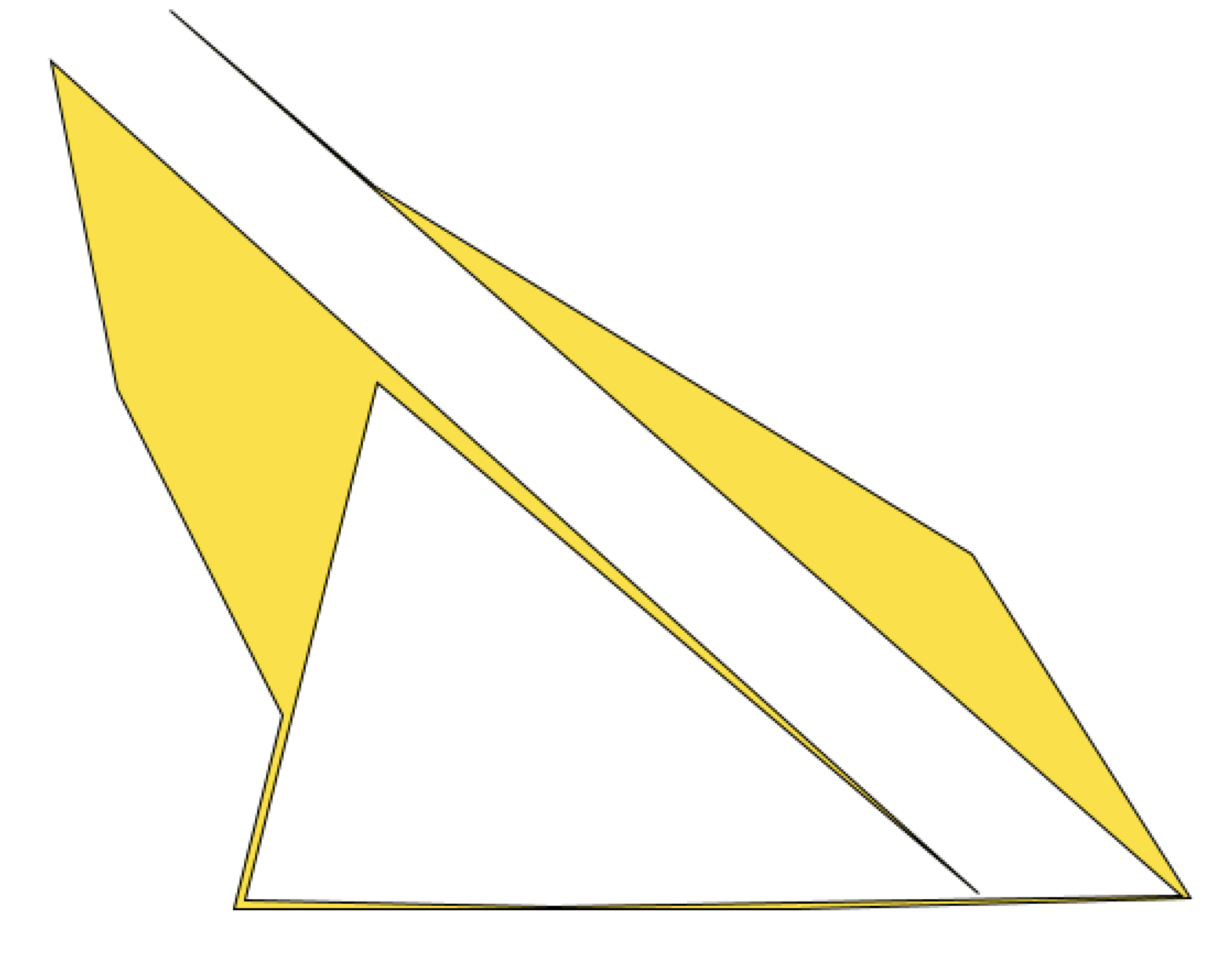}}
\resizebox{!}{2.4in}{\includegraphics{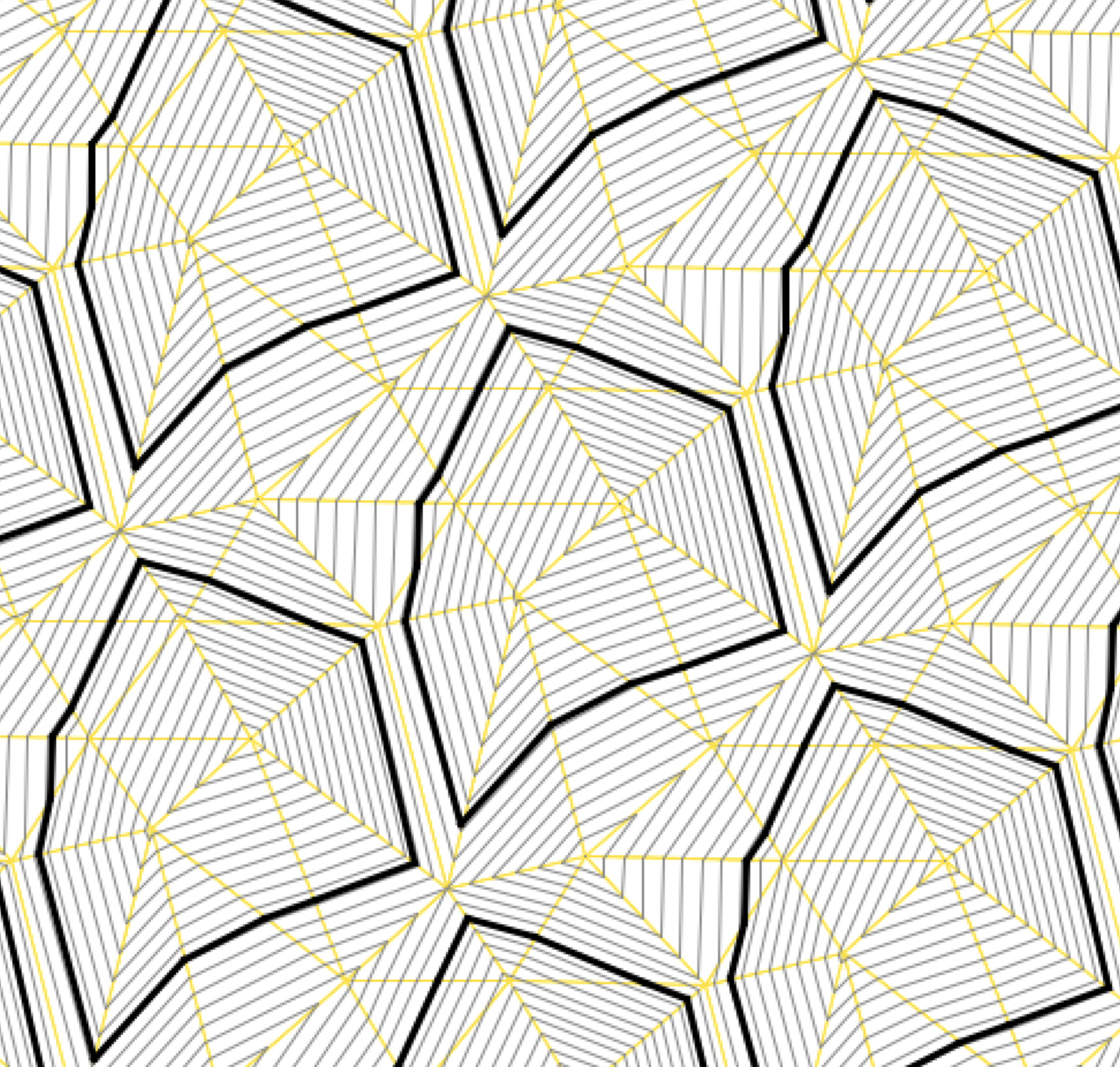}}
\newline
{\bf Figure 4.7:\/} Extrinsic and Intrinsic structures compared
\end{center}

Once we fix this direction we get a {\it level set foliation\/}
of the pup tent
by curves that lie in planes perpendicular to this direction.
On the right side of Figure 4.7, we have drawn the level set foliation in the
intrinsic flat structure and then
lifted the picture to the universal cover.
The highlighted loop on the right
corresponds to the loop shown on the left.

My supplementary notes [{\bf S\/}]
give additional information about the pup tent, including
a pattern than you can cut out and (with effort) fold into one.
Also, my Java program lets you see many more pictures like this.

\newpage

  \section{References}

\setlength{\parskip}{1\baselineskip}

\noindent
[{\bf ALM\/}] P. Arnoux, S. Leli\`evre, A. M\'alaga, {\it Diplotori:
  a family of polyhedral flat tori.\/} in preparation

\noindent
[{\bf Br\/}], U. Brehm, Oberwolfach report (1978)

\noindent
[{\bf BE\/}] J. Bokowski and A. Eggert,
{\it All realizations of M\"obius' Torus with
  7 Vertices\/}, topologie structurale n\'um 17, (1991)

\noindent
[{\bf BZ1\/}] Y. D. Burago and V. A. Zalgaller,
{\it Polyhedral realizations of developments\/} (in Russian)
Vestnik Leningrad Univ. 15, pp 66--80 (1960)

\noindent
[{\bf BZ2\/}] Y. D. Burago and V. A. Zalgaller,
{\it Isometric Embeddings of Two Dimensional Manifolds with a
  polyhedral metric into $\R^3$\/}, Algebra i analiz 7(3) pp 76-95
(1995)  Translation in St. Petersburg Math Journal (3)3, pp 369--385

\noindent
[{\bf Cs\/}]
Á. Császár, ``A polyhedron without diagonals,'' 
\textit{Acta Sci. Math. (Szeged)}, vol.\ 13, pp.\ 140--142, 1949.

\noindent
[{\bf DS\/}]
P. Doyle and R. E. Schwartz
{\it Collapsibility and Near Universality for Vertex-Minimal Paper
  Tori\/},
arXiv 2510.112623 (2025)

\noindent
[{\bf E\/}]
M. Ellison, ``Realizing abstract simplicial complexes with specified
edge lengths'' arXiv 2312.05376 (2023)

\noindent
[{\bf G\/}], M. Gardner, {\it Mathematical Games: On the remarkable
 Cs\'asz\'ar  polyhedron and its applications in problem solving\/},
Scientific American {\bf 232\/}, 5 (1975) pp 102-107.

\noindent
[{\bf HLZ\/}], S. Hougardy, F H. Lutz, and M. Zelke,
{\it Polyhedral Tori with Minimal Coordinates\/},  arXiv 0709.2794

\noindent
[{\bf LT\/}] F. Lazarus, F. Tallerie, {\it A Universal Triangulation
  for Flat Tori\/}, CS arXiv 2203.05496 (2024)

\noindent
[{\bf L\/}] F. H. Lutz, {\it Cs\'asz\'ar's Torus\/}, Electronic
  Geometry models (2002) 2001.02.069
  
  \noindent
  [{\bf QC\/}], P. T. Quintanar Cort\'es,
  {\it Plongements poly\'edriques du tore carr\'e plat\/},
  PhD thesis, Universit\'e Claude Bernard Lyon~1, 2019.
\newline  
  Available at {\tt http://www.theses.fr/2019LYSE1354}.

\noindent
[{\bf S\/}] R. E. Schwartz,  {\it Notes on the Pup Tent\/}, informal
supplementary  notes
(2025) \newline
Available at {\tt http://www.math.brown.edu/$\sim$res/MathNotes/puptent.pdf\/}

\noindent
[{\bf Se\/}] H. Segerman, {\it Visualizing Mathematics with 3D Printing\/},
Johns Hopkins U. Press  (2016)

\noindent
[{\bf T\/}] T. Tsuboi, {\it On Origami embeddings of flat tori\/},
arXiv 2007.03434 (2020)
  
\noindent
[{\bf Tu\/}] V. Tugay\'e, personal communication (2025)

\noindent
[{\bf Z\/}] Z. Zou, {\it Vertex Minimal Hyperbolic Origami 2-Torus\/}, arXiv
2509.18668 (2025)

\end{document}